\documentclass[a4paper,12pt]{article}
\usepackage{latexsym}
\newtheorem{theorem}{Theorem}[section]
\newtheorem{corollary}[theorem]{Corollary}

\newtheorem{example}[theorem]{Example}
\newtheorem{proposition}[theorem]{Proposition}
\newtheorem{remark}[theorem]{Remark}

\title{On locally $LC$-spaces\thanks{1991 Math.\ Subject
Classification --- Primary: 54D20, 54A05, Secondary: 54D45,
54G99.\protect\newline Key words and phrases --- $LC$-space,
locally $LC$-space, Lindel\"{o}f set.}}
\author{Julian {\sc Dontchev}, Maximilian {\sc Ganster} and Alev
{\sc Kanibir}}
\date{}
\begin{document}
\baselineskip=20pt plus 1pt minus 1pt
\newcommand{\fxy}{$f \colon (X,\tau) \rightarrow (Y,\sigma)$}
\newcommand{\llc}{locally $LC$-space}
\maketitle
\begin{abstract}
A topological space $(X,\tau)$ is called a \llc\ if every point
of $X$ has a neighborhood $U$ such that every Lindel\"{o}f subset
of $(U,\tau|U)$ is a closed subset of $(U,\tau|U)$. The aim of
this paper is to continue the study of \llc s.
\end{abstract}

\section{Introduction}\label{s1}

Classical generalizations of Lindel\"{o}f spaces such as
hereditarily Lindel\"{o}f and maximal Lindel\"{o}f spaces have
had their major impact in the development of General Topology.
A certain class of spaces, relatively new as a concept but
extensively studied in recent years, is the class of $LC$-spaces.
A topological space $(X,\tau)$ whose Lindel\"{o}f subsets are
closed is called an {\em $LC$-space} by Gauld, Mr\v{s}evi\'{c},
Reilly and Vamanamurthy \cite{GMRV1} and by Mukherji and Sarkar
\cite{MS1}. This concept emerged from the study of maximal
Lindel\"{o}f spaces \cite{C1} as being a notion having a close
relationship to P-spaces.

$LC$-spaces generalize Wilansky's KC-spaces \cite{W1} and
Hausdorff P-spaces \cite{M1}. On the other hand every $LC$-space
is a cid-space, i.e., all countable sets are closed and discrete,
and hence $T_1$ and anti-compact (= pseudo-finite) \cite{B1}. An
extensive list of references on $LC$-spaces as well as some
generalizations of the concept can be found in \cite{DGK1}.

Recently, Ganster, Kanibir and Reilly \cite{GKR1} introduced the
class of \llc s. By definition, a topological space $(X,\tau)$
is called a {\em locally $LC$-space} if each point of $X$ has a
neighborhood which is an $LC$-subspace. In \cite{GKR1}, the
authors proved that a space $(X,\tau)$ is an $LC$-space if and
only if each point of $X$ has a closed neighborhood that is an
$LC$-subspace. Thus every regular \llc\ is an $LC$-space, a
result first proved by Hdeib and Pareek in \cite{HP1}. The
following example shows that we cannot replace `regular' by
`Hausdorff'.

\begin{example}\label{e1}
{\em \cite{GJ1,GKR1} There exists a Hausdorff, \llc\ $(X,\tau)$,
which is not an $LC$-space. Let $Z$ be a set of cardinality
$\aleph_{1}$ with a distinguished point $z_0$. The topology on
$Z$ is defined as follows: each $z \not= z_0$ is isolated while
the basic neighborhoods of $z_0$ are the cocountable subsets of
$Z$ containing $z_0$. Note that $Z$ is a Lindel\"{o}f
$LC$-space. The space $(X,\tau)$ will be constructed from copies
of $Z$. For each $n \in \omega$, let $X_n$ be a copy of $Z$,
where $x_n$ denotes the non-isolated point of $X_n$. Let $X^{*}
= \sum_{n \in \omega} X_n$ denote the topological sum of the
spaces $X_n$ and let $X = X^{*} \cup \{p\}$ with $p \not\in
X^{*}$. A topology $\tau$ on $X$ can be defined if, in addition,
we specify the basic open neighborhoods of $p$. They are the
union of $\{p\}$ and a cocountable subset of $\cup \{ X_n
\setminus \{ x_n \} \colon n \geq k \}$ for some $k \in \omega$.
$(X,\tau)$ is a Hausdorff space that fails to be an $LC$-space
\cite{GJ1}. However, as shown in \cite{GKR1}, $(X,\tau)$ is a
\llc.}
\end{example}

\section{Locally $LC$-spaces}\label{s2}

\begin{proposition}\label{t1}
For a topological space $(X,\tau)$ the following conditions are
equivalent:

{\rm (1)} $X$ is a \llc s.

{\rm (2)} Every point of $X$ has an open neighborhood, which is
an $LC$-subspace of $X$.
\end{proposition}

{\em Proof.} Follows from the fact that every subspace of an
$LC$-space is an $LC$-space. $\Box$

\begin{proposition}\label{t2}
Every subspace of a \llc\ is a \llc.
\end{proposition}

{\em Proof.} Let $(X,\tau)$ be a \llc\ and let $A \subseteq X$.
By assumption, for each $x \in A$, there exists $U \in \tau$ such
that $(U,\tau|U)$ is an $LC$-space. Note that $V = U \cap A$ is
an open neighborhood of $x$ in $(A,\tau|A)$ and that $(V,\tau|V)$
is an $LC$-subspace. By Proposition~\ref{t1}, $(A,\tau|A)$ is a
\llc. $\Box$

\begin{proposition}\label{t3}
If $(X,\tau)$ has an open cover by locally LC-subspaces, then $X$
is a \llc.
\end{proposition}

{\em Proof.} Let  $X = \cup_{i \in I} G_i$ be on open cover of
$X$ where each $G_i$ is a \llc, and let $x \in X$. Choose $j \in
I$ such that $x \in G_j$. If $U_j$ is an open neighborhood of $x$
in $G_j$ such that $U_j$ is an $LC$-subspace of $G_j$,
then $U_j$ is also open in $(X,\tau)$. By Proposition~\ref{t1},
$(X,\tau)$ is a \llc. $\Box$

\begin{corollary}\label{c1}
Let $(X_{\alpha},{\tau}_{\alpha})_{\alpha \in \Omega}$ be a
family of topological spaces. For the topological sum $X =
\sum_{\alpha \in \Omega} X_{\alpha}$ the following conditions are
equivalent:

{\rm (1)} $X$ is a \llc.

{\rm (2)} Each $X_{\alpha}$ is a \llc.
\end{corollary}

{\em Proof.} Follows from Proposition~\ref{t2} and
Proposition~\ref{t3}.
$\Box$

\begin{proposition}\label{tpr1}
Let $(X_{i},\tau_{i})_{i \in F}$ be a finite family of Hausdorff
spaces. If each $X_i$ is a \llc, then the product space $X =
\prod_{i \in F} X_i$ is also a \llc.
\end{proposition}

{\em Proof.} Let $(X,\tau)$ and $(Y,\sigma)$ be Hausdorff
\llc s. Let $(x,y) \in X \times Y$. By assumption, there exists
$U \in \tau$ and $V \in \sigma$ such that $x \in U$, $y \in V$
and both $U$ and $V$ are $LC$-subspaces of $X$ and $Y$,
respectively. By \cite[Theorem 2]{DG1}, $U \times V$ is an
$LC$-subspace of the product space $X \times Y$. By
Proposition~\ref{t1}, $X \times Y$ is a \llc. $\Box$

\begin{remark}\label{r1}
{\em We note that the Hausdorff condition can be reduced to the
weaker separation property $R_1$. Recall that a space
$(X,\tau)$ is called an {\em $R_1$-space} if $x$ and $y$ have
disjoint neighborhoods whenever ${\rm cl} \{x\} \not= {\rm cl}
\{y\}$. Clearly, a space is Hausdorff if and only if it is $T_1$
and $R_1$.}
\end{remark}

{\bf Question.} Does Proposition~\ref{tpr1} remain true if we
drop the requirement that the spaces in question have to be
Hausdorff (or $R_1$)?

\begin{proposition}\label{t31}
Every locally hereditarily Lindel\"{o}f, \llc\ is discrete.
\end{proposition}

{\em Proof.} Let $(X,\tau)$ be a locally hereditarily
Lindel\"{o}f and a locally $LC$-space. We may assume that every
point $x \in X$ has an open neighbourhood $W$ that is both
hereditarily Lindel\"{o}f and an $LC$-space. But this means that
$W$ is an open discrete subspace of $(X,\tau)$. Thus $(X,\tau)$
is discrete as well. $\Box$

\begin{corollary}
A \llc\ is discrete if and only if it is locally finite. $\Box$
\end{corollary}

\begin{proposition}\label{ppn1}
Every \llc\ is a $T_1$-space.
\end{proposition}

{\em Proof.} Assume that $(X,\tau)$ is a \llc\ and let $x \in X$.
For every $y \not= x$, there exists $U \in \tau$ such that $y \in
U$ and $(U,\tau|U)$ is an $LC$-space and hence a $T_1$-space.
Clearly, $y \in U \setminus \{ x \} \in \tau$ and $x \not\in U$.
This shows that $X$ is $T_1$. $\Box$

\begin{proposition}\label{ppn2}
If every $LC$-subspace of every Lindel\"of subset of a
topological space $(X,\tau)$ is Lindel\"of, then $X$ is a \llc\
if and only if $X$ is an $LC$-space.
\end{proposition}

{\em Proof.} Assume that $X$ is a \llc. Let $A \subseteq X$ be
Lindel\"of and let $x \not\in A$. Since $X$ is an \llc, there
exists $U \in \tau$ such that $x \in U$ and $(U,\tau|U)$ is an
$LC$-space. Since every subspace of an $LC$-space is an
$LC$-space, $U \cap A$ is an $LC$-space. By assumption, $U \cap
A$ is Lindel\"of and hence closed in $(U,\tau|U)$. Thus, $U
\setminus A$ is open in $(X,\tau)$, contains $x$ and is disjoint
from $A$. This shows that $A$ is closed and consequently $X$ is
an $LC$-space. $\Box$

\begin{proposition}\label{t4}
Open (and hence also closed) bijective images of \llc s are \llc
s.
\end{proposition}

{\em Proof}. Let \fxy\ be open and bijective and let $(X,\tau)$
be a \llc. Let $y \in Y$. Choose $x \in X$ such that $f(x) = y$.
Since $(X,\tau)$ is a \llc, then there exists $U \in \tau$ such
that $x \in U$ and $(U,\tau|U)$ is an $LC$-space. Since $f$ is
open, then $f(U)$ is an open neighborhood of $y$ in $(Y,\sigma)$.
Since open, bijective images of $LC$-spaces are also $LC$-spaces,
then $(f(U),\sigma|{f(U)})$ is an $LC$-subspace of $(Y,\sigma)$.
By Proposition~\ref{t1}, $(Y,\sigma)$ is a \llc. $\Box$

\begin{corollary}\label{c2}
The property `locally $LC$-space' is a topological property.
$\Box$
\end{corollary}

\
\begin{center}
Department of Mathematics\\PL 4, Yliopistonkatu 5\\
University of Helsinki\\00014 Helsinki\\Finland\\e-mail:
{\tt dontchev@cc.helsinki.fi}
\end{center}
\
\begin{center}
Department of Mathematics\\Graz University of
Technology\\Steyrergasse 30\\A-8010 Graz\\Austria\\e-mail:
{\tt ganster@weyl.math.tu-graz.ac.at}
\end{center}
\
\begin{center}
Department of Mathematics\\Faculty of Science\\Haceteppe
University\\06532 Beytepe/Ankara\\Turkey\\e-mail: {\tt
kanibir@eti.cc.hun.edu.tr}
\end{center}
\
\end{document}